\documentclass[11pt,fleqn]{article}
\usepackage{amsfonts, epsfig}
\newcommand{\ol}{\setlength{\itemsep}{0pt.}\begin{enumerate}}
\newcommand{\eol}{\end{enumerate}\setlength{\itemsep}{-\parsep}}
\newcommand{\ignore}[1]{}
\newcommand{\eqdef}{\stackrel{\rm def}{=}}
\title{Linear programming bounds for codes via a covering argument}
\author{Michael Navon\thanks{School of Computer Science and
    Engineering, Hebrew University, Jerusalem, Israel.}
\and Alex Samorodnitsky\thanks{School of Computer Science and
    Engineering, Hebrew University, Jerusalem, Israel.}}
\begin{document}
\date{}
\maketitle
 
 
\newtheorem{THEOREM}{Theorem}[section]
\newenvironment{theorem}{\begin{THEOREM} \hspace{-.85em} {\bf :} 
}%
                        {\end{THEOREM}}
\newtheorem{LEMMA}[THEOREM]{Lemma}
\newenvironment{lemma}{\begin{LEMMA} \hspace{-.85em} {\bf :} }%
                      {\end{LEMMA}}
\newtheorem{COROLLARY}[THEOREM]{Corollary}
\newenvironment{corollary}{\begin{COROLLARY} \hspace{-.85em} {\bf 
:} }%
                          {\end{COROLLARY}}
\newtheorem{PROPOSITION}[THEOREM]{Proposition}
\newenvironment{proposition}{\begin{PROPOSITION} \hspace{-.85em} 
{\bf :} }%
                            {\end{PROPOSITION}}
\newtheorem{DEFINITION}[THEOREM]{Definition}
\newenvironment{definition}{\begin{DEFINITION} \hspace{-.85em} {\bf 
:} \rm}%
                            {\end{DEFINITION}}
\newtheorem{EXAMPLE}[THEOREM]{Example}
\newenvironment{example}{\begin{EXAMPLE} \hspace{-.85em} {\bf :} 
\rm}%
                            {\end{EXAMPLE}}
\newtheorem{CONJECTURE}[THEOREM]{Conjecture}
\newenvironment{conjecture}{\begin{CONJECTURE} \hspace{-.85em} 
{\bf :} \rm}%
                            {\end{CONJECTURE}}
\newtheorem{MAINCONJECTURE}[THEOREM]{Main Conjecture}
\newenvironment{mainconjecture}{\begin{MAINCONJECTURE} \hspace{-.85em} 
{\bf :} \rm}%
                            {\end{MAINCONJECTURE}}
\newtheorem{PROBLEM}[THEOREM]{Problem}
\newenvironment{problem}{\begin{PROBLEM} \hspace{-.85em} {\bf :} 
\rm}%
                            {\end{PROBLEM}}
\newtheorem{QUESTION}[THEOREM]{Question}
\newenvironment{question}{\begin{QUESTION} \hspace{-.85em} {\bf :} 
\rm}%
                            {\end{QUESTION}}
\newtheorem{REMARK}[THEOREM]{Remark}
\newenvironment{remark}{\begin{REMARK} \hspace{-.85em} {\bf :} 
\rm}%
                            {\end{REMARK}}
 
\newcommand{\thm}{\begin{theorem}}
\newcommand{\lem}{\begin{lemma}}
\newcommand{\pro}{\begin{proposition}}
\newcommand{\dfn}{\begin{definition}}
\newcommand{\rem}{\begin{remark}}
\newcommand{\xam}{\begin{example}}
\newcommand{\cnj}{\begin{conjecture}}
\newcommand{\mcnj}{\begin{mainconjecture}}
\newcommand{\prb}{\begin{problem}}
\newcommand{\que}{\begin{question}}
\newcommand{\cor}{\begin{corollary}}
\newcommand{\prf}{\noindent{\bf Proof:} }
\newcommand{\ethm}{\end{theorem}}
\newcommand{\elem}{\end{lemma}}
\newcommand{\epro}{\end{proposition}}
\newcommand{\edfn}{\bbox\end{definition}}
\newcommand{\erem}{\bbox\end{remark}}
\newcommand{\exam}{\bbox\end{example}}
\newcommand{\ecnj}{\bbox\end{conjecture}}
\newcommand{\emcnj}{\bbox\end{mainconjecture}}
\newcommand{\eprb}{\bbox\end{problem}}
\newcommand{\eque}{\bbox\end{question}}
\newcommand{\ecor}{\end{corollary}}
\newcommand{\eprf}{\bbox}
\newcommand{\beqn}{\begin{equation}}
\newcommand{\eeqn}{\end{equation}}
\newcommand{\wbox}{\mbox{$\sqcap$\llap{$\sqcup$}}}
\newcommand{\bbox}{\vrule height7pt width4pt depth1pt}
\newcommand{\qed}{\bbox}
\def\sup{^}

\def\H{\{0,1\}^n}

\def\S{S(n,w)}

\def \E{\mathbb E}
\def \R{\mathbb R}
\def \Z{\mathbb Z}
\def \F{\mathbb F}

\def\<{\left<}
\def\>{\right>}
\def \({\left(}
\def \){\right)}
\def \e{\epsilon}
\def \r{\rfloor}
\def \S{{\cal S}}
\def \b{\bullet}

\def \L{\Lambda}
\def \G{{\Gamma}}

\def \O{\Omega}
\def \D{\Delta}
\def \d{\delta_0}
\def \l{\lambda}

\def \grad{\bigtriangledown}

\newcommand{\rarrow}{\rightarrow}
\newcommand{\lrarrow}{\leftrightarrow}

\begin{abstract}
We recover the first linear programming bound of McEliece, Rodemich,
Rumsey, and Welch for binary
error-correcting codes and designs via a covering argument. It is possible 
to show, interpreting the following notions appropriately, that if a
code has a large distance, then its dual has a small covering radius
and, therefore, is large. This implies the original code to be small.

We also point out (in conjunction with further work) that this bound
is a natural isoperimetric constant of the Hamming cube, related to
its Faber-Krahn minima.

While our approach belongs to the general framework of Delsarte's linear
programming method, its main technical ingredient is
Fourier duality for the Hamming cube. In particular,
we do not deal directly with Delsarte's linear
program or orthogonal polynomial theory.

\end{abstract}

\section{Introduction}
This paper takes another look at the first linear programming bound
on binary error correcting codes, or, alternatively, on 
optimal packing of Hamming balls in a Hamming cube. 

The bound was originally proved by McEliece, Rodemich, Rumsey, and
Welch \cite{mrrw}, following Delsarte's linear programming approach
\cite{dels}. 
Delsarte showed the distance distribution of a binary code 
to satisfy a family of linear constraints whose coefficients can be
viewed as values of a certain family of orthogonal polynomials, i.e.,
the Krawchouk polynomials. This made it possible to construct a linear
programming relaxation 
of the original combinatorial question, and to view the obtained
linear program as an extremal problem 
in orthogonal polynomials. Good feasible solutions of the {\it dual
  program} were  constructed in \cite{mrrw} using tools from the theory of
orthogonal polynomials. These solutions lead to the bound, known
as {\it the first linear programming bound} (or the first JPL
bound). This bound is the  best known upper bound on cardinality of a
code with a given minimal distance, for a significant range of distances.

Delsarte's approach extends to a family of finite metric spaces, known
as {\it commutative association schemes}. A Hamming cube is one example
of an association scheme. Another relevant example is the Hamming
sphere. In \cite{mrrw} 
good feasible solutions to Delsarte's linear program for the Hamming
sphere are constructed. These lead to best known upper bounds on
constant weight error correcting codes (ball packing in the Hamming
sphere), and, combined with the Bassalygo-Elias inequality, to the
best known upper bound on binary codes. This bound is known as the
{\it second linear programming bound}.

We refer to \cite{mrrw, lev-survey, BI, MS} for a detailed exposition
of the notions discussed above, including error-correcting codes and
their significance, packing in metric spaces, association schemes,
Delsarte's linear program, and orthogonal polynomials.

The point of view presented in this paper is somewhat different. Our
main tool is Fourier 
analysis on the group $\F^n_2$, or, equivalently, on the Hamming cube
$\H$. We follow the approach of Kalai and Linial \cite{kl} in which
the characteristic 
function of a binary code is viewed as a real-valued function on the
cube. A study of the Fourier transform of this function and its simple 
by-products makes it possible to recover Delsarte's linear program in
a form which does not require treatment of Krawchouk polynomials. 

Moreover, this viewpoint allows an easy access to new geometric
information. Specifically, we establish a simple relation between the
minimal distance (equivalently, packing 
radius) of a code and the {\it essential covering radius} of its
dual. Recall that $r$ 
is a covering radius of a subset $C$ of $\H$ if
the union of Hamming 
balls of radius $r$ centered at the points of $C$ covers the whole
space. Here we use a somewhat weaker notion, and require this union
of balls to cover only a significant fraction of the space. 

This observation, which we consider to be the main contribution of this
paper, leads to a simple proof of the first linear programming
bound. In particular, we do not need to deal directly with Delsarte's
linear program or orthogonal polynomial theory.

We move to the principal definitions and to the statement of the main
results.

{\it A binary error-correcting code with block length $n$ and 
minimal distance $d$} is a subset of the $n$-dimensional Hamming cube
in which the distance between any two distinct points is at least
$d$.
Let $A(n,d)$ be the maximal size of such a code. 
In this paper we are interested in the case in which the distance
$d$ is linear in the length $n$ of the code, and we let the length $n$
go to infinity. In this case $A(n,d)$ is known \cite{lev-survey} to
grow exponentially in $n$, and we consider the quantity 
$$
R(\delta) = \limsup_{n\rightarrow \infty} \frac 1n
\log_2 A(n,\lfloor \delta n \rfloor),
$$
also known as {\it the asymptotic maximal rate of the code
with relative distance $\delta$} for $0 \le \delta \le \frac12$.

Next, we need the notion of a {\it maximal eigenvalue} of a subset of
the cube. We say that two elements $x, y$ of $\H$ are adjacent and
write $x \sim y$ if the Hamming distance between $x$ and $y$ is $1$. 
Let $A$ be the adjacency matrix of the obtained graph. For $B
\subseteq \H$, set  
$$
\l_B = \max\left\{\frac{\<Af,f\>}{\<f,f\>};~~f:~\H \rarrow
\R,~\mbox{supp}(f) \subseteq B\right\}
$$ 
In other words, $\l_B$ is the maximal eigenvalue of adjacency matrix
of the subgraph of $\H$ induced by the vertices of $B$. 

Our main technical claim is 
\pro
\label{pro:ineq}
Let $C$ be a code with block length $n$ and minimal distance $d$. Let
$B$ be a subset of $\H$ with $\l_B \ge n - 2d + 1$. Then
$$
|C| \le n |B|
$$
\epro

A linear code $C$ is a linear subspace of $\F^n_2$. The dual
code $C^{\perp}$ is the orthogonal subspace, that is it contains all
the vectors orthogonal to $C$ over $\F_2$. Proposition~\ref{pro:ineq}
has an appealing geometric interpretation for linear codes.  
\pro
\label{pro:cover}
Let $C$ be a linear code with block length $n$ and minimal distance $d$. Let
$B$ be a subset of $\H$ with $\l_B \ge n - 2d + 1$. Then 
$$
\Big | \bigcup_{z \in C^{\perp}} \(z + B\) \Big | \ge \frac{2^n}{n}
$$
\epro

In other words, replacing every point in the dual code by a
(shifted) copy of $B$, we will cover a large fraction of the space
$\H$. Proposition~\ref{pro:ineq} for linear codes is an immediate
corollary of (\ref{cover}) since $|C|\cdot |C^{\perp}| = 2^n$.
  
A code $C'$ has {\it dual distance} $d$ if Fourier transform of its
characteristic function vanishes on points of Hamming weight $0 < |S|
< d$. In particular, the dual distance of a linear code is easily seen
to equal the minimal
distance of its dual (cf. discussion in
Subsection~\ref{subsec:fourier}). Hence, the following claim
generalizes 
Proposition~\ref{pro:cover}. 
\pro
\label{pro:dual}
Let $C'$ be a code with block length $n$ and dual distance $d$. Let
$B$ be a subset of $\H$ with $\l_B \ge n - 2d + 1$. Then 
\beqn
\label{cover}
\Big | \bigcup_{z \in C'} \(z + B\) \Big | \ge \frac{2^n}{n}
\eeqn
\epro

Hamming balls are a good choice for the covering set $B$. 
\lem
\label{lem:eigen-ball}
Let $B(r)$ be a Hamming ball of radius $r$. Then 
$$
\l_{B(r)} \ge 2\sqrt{r(n-r)} - o(n) 
$$
\elem

Proposition~\ref{pro:dual} together with Lemma~\ref{lem:eigen-ball}
lead to a relation between the dual distance of a
code and its essential covering radius.
\cor
\label{cor:cover}
Let $C'$ be a code with block length $n$ and dual distance $d$. Then
the essential covering radius of $C'$ is at most
\beqn
\label{cover-rad}
r \le \frac n2 - \sqrt{d(n-d)} + o(n) 
\eeqn
In particular, 
let $C$ be a linear code with block length $n$ and minimal distance $d$. Then
the essential covering radius of the dual code $C^{\perp}$ is at most
$
r \le \frac n2 - \sqrt{d(n-d)} + o(n). 
$
\ecor

Proposition~\ref{pro:ineq} together with Lemma~\ref{lem:eigen-ball}
lead to an upper bound on the size of a code $C$ with block
length $n$ and minimal distance $d$. They show that there
exists a radius $r \le \frac n2 - \sqrt{d(n-d)} + o(n)$ such that 
\beqn
\label{ineq:firstLP}
|C| \le n |B(r)| 
\eeqn
Corollary~\ref{cor:cover} gives a geometric explanation of this bound
for a linear code $C$.  The balls of radius $r$
centered at the points of the dual code $C^{\perp}$ cover an
$\(1/n\)$-fraction of the space. Therefore $| C^{\perp}| \cdot
|B(r)| \ge 2^n/n$, and $|C| = 2^n/|C^{\perp}| \le n
|B(r)|$. This allows us to view the bound (\ref{ineq:firstLP}) as a
{\it covering bound}. 

For a general code the covering interpretation of (\ref{ineq:firstLP})
is more tenuous since, in particular, there is no natural notion of
the dual code. However, the analytic reasoning leading to
(\ref{ineq:firstLP}) can be viewed as a functional version of the
covering argument above (see Subsection~\ref{subsec:gencod}).

The cardinality of a Hamming ball of radius
$r$ is $2^{n\(H\(r/n\) + o(1)\)}$ \cite{vL}. Substituting the value $r
= \frac n2 - \sqrt{d(n-d)} + o(n)$ on the right hand side of
(\ref{ineq:firstLP}), we have
$$
|C| \le 2^{n\(H\(1/2 - \sqrt{d/n\(1-d/n\)}\)   + o(1)\)}
$$
This bounds the asymptotic maximal rate of a
code with relative distance $\delta$, 
$$
R(\delta) \le  H\(1/2 - \sqrt{\delta\(1-\delta\)}\)
$$
This is the  first linear programming bound for error-correcting
codes. 

Finally, a code with dual distance $d$ is 
a {\it design of strength $d$}
\cite{lev-survey} ~(or a {\it $(d-1)$-wise independent set} \cite{abi}).
Proposition~\ref{pro:dual} together with Lemma~\ref{lem:eigen-ball}
lead to the {\it first linear programming bound} for designs \cite{mrrw}. 

\noindent {\bf Summing up} 

Three notions of duality are relevant to this discussion. The first is
linear programming duality as represented by the 
primal and dual linear programs of Delsarte. Recall that the primal linear
program of Delsarte is a relaxation of the combinatorial question on
cardinality of an optimal code. The linear programming bounds on
codes are obtained by constructing good feasible solutions for the
dual program. 

The second notion is the Fourier duality, illustrated by
the Kalai-Linial approach to the problem. Viewing the
characteristic function of a code as a real-valued function on the
cube, and studying the properties of this function and its Fourier
transform lead to an equivalent version of Delsarte's linear program. 

The third notion is the duality between packing and covering problems
in hypergraphs \cite{ael}. The vertices of the pertinent hypergraph
are the vertices of the cube and the edges are Hamming balls. The fractional
packing and covering problems are dual linear programs. This induces a
duality relation between their integer versions which
are of interest here. Generally, covering is much easier than
packing. For instance, integrality gap for covering is logarithmic at worst
\cite{lov}, while for packing it could be much larger \cite{ael}. In
the context of coding theory, the asymptotics of optimal packings are
unknown, while the asymptotics of optimal coverings are easy
to find \cite{CHLL}.

The main observation of this paper is that, in our case, Fourier
duality makes it 
possible to pass from a ``hard' packing problem of finding the
maximal cardinality of a code with a given minimal distance to an
``easy'' covering question of determining the minimal cardinality of a
code with a given covering radius. We suggest that this point of view
might explain the power of the resulting bound, namely the first linear
programming bound for error-correcting codes.  

We also point out that this bound is a natural isoperimetric constant of 
the Hamming cube, related to its {\it Faber-Krahn minima} (\cite{ft, sam:sob},
see the discussion below).  

\noindent {\bf Related work}
\begin{enumerate}
\item
Our research was motivated by a recent result of 
Friedman and Tillich \cite{ft}. Using 
methods from algebraic graph theory the authors prove the first linear
programming bound for linear binary codes. In particular,
Proposition~\ref{pro:ineq} for linear codes and
Lemma~\ref{lem:eigen-ball} are proved in \cite{ft}.\footnote{We give a
different proof of Lemma~\ref{lem:eigen-ball}, based on an explicit
construction of a function with a large Rayleigh ratio.} The appeal of
\cite{ft} is in suggesting a way to work with Delsarte's linear
inequalities without resorting to the language and tools of orthogonal
polynomial theory. 
\item
Combining the approach of Friedman and Tillich with the
Fourier-analytic view of  Delsarte's linear program due to Kalai and
Linial allowed us to extend this approach to general binary codes,
with, we believe, a simpler proof. After completing our work on the
conference version of this paper \cite{ns}, we learned that Fourier
analysis was used in a similar manner by Cohn and Elkies
\cite{ce} to give a simpler proof of Levenshtein's bound on sphere
packing in $\R^n$. In particular, \cite{ce} contains (somewhat
implicitly) arguments analogous to our proofs of
Proposition~\ref{pro:ineq} and Lemma~\ref{lem:eigen-ball}.
\item
The relation between the dual distance of a code and its covering
radius has been extensively investigated in the coding literature
(see \cite{tie, ahll} and the references there). The best known bounds
are obtained via linear
programming approach and are somewhat weaker than
(\ref{cover-rad}). This, of course, stands to reason, since
covering radius of a code is, in general, larger that the essential covering
radius. The best known upper bound on the covering
radius $r_c$ of a code with dual distance $d$ is \cite{tie} 
$$
r_c \le \frac n2 -  \sqrt{\frac{d}{2}\(n-\frac{d}{2}\)} + o(n)
$$
Better bounds are known for linear codes \cite{ahll}.
\end{enumerate}

\noindent {\bf Extensions and ramifications}
\label{subsec:extend}
\begin{enumerate}
\item
The approach of this paper extends to general
commutative association schemes \cite{sam:sch}. Given an association
scheme with $k+1$ classes, it is possible to a formal ``Fourier transform'' on
vectors in $\R^{k+1}$. This transform,  contrary to a genuine
Fourier transform, is not self-dual. The solution is to
define a pair of 'inverse transforms'. An appropriate pair of linear
transformations is given by the transition matrices of
the scheme. It is possible to define a 
convolution operation for each of these transforms, which is
commutative and associative, and which is taken by the transform to a
point-wise multiplication. This allows to recover the best
known bounds for codes and designs in commutative association schemes
via a Fourier-analytic proof similar to the proof of
Proposition~\ref{pro:ineq}. 
\item
Friedman and Tillich \cite{ft} ask what are the optimal covering
subsets of the Hamming cube, in the sense of
Proposition~\ref{pro:ineq}. This
question is answered in \cite{sam:sob}, by way of a modified
logarithmic Sobolev inequality for the Hamming cube. If $B$ is a
Hamming ball and $X$ is a 
subset of $H$ with $|X| = |B|$, then $\l_X \le (1+o(1)) \cdot
\l_B$. Hence Hamming balls are asymptotically optimal. In the
terminology of \cite{ft}, Hamming balls are the Faber-Krahn minimizers
for the Hamming cube (up to a negligible error).   
\end{enumerate}
\section{The proofs}
\subsection{Fourier analysis on $\F^n_2$}
\label{subsec:fourier}
We refer to \cite{kkl} for background in Fourier analysis on
$\F^n_2$. Here we list several necessary definitions and simple
facts. 

$\F_2^n$ is a finite Abelian group, therefore its characters
$\{W_S\}_{S \in \F_2^n}$ constitute a group (the
{\em dual group} which is isomorphic to $\F_2^n$.)
The character $W_S$ is
a function from $\F_2^n$ to $\{-1,1\}$,
defined as: $W_S(x) = (-1)^{\<x,S\>} $.
The characters $\{W_S\}_{S \in \F_2^n}$ form an orthonormal basis 
in the space of real-valued functions on $\F_2^n$, equipped with
uniform probability distribution.  

Write $\E f$ for $\frac{1}{2^n} \sum_{x \in \F^n_2} f(x)$. 
For $f: \F_2^n \rightarrow \R$, define 
$\widehat f: \F_2^n \rightarrow \R$ as  $\widehat f(S) = \<f,W_S\> \eqdef 
\E \(f \cdot W_S\)$.
The function $\widehat f$ is the {\em Fourier Transform} of $f$. The
Parseval identity states $\E f g = \<f,g\> = \<\widehat{f},\widehat{g}\>
\eqdef \sum \widehat{f} \widehat{g}$.

For $f, g: \F_2^n \rightarrow \R$, the {\it convolution} of $f$ and $g$
is defined by $(f \ast g)(x) = \E_{y} f(y)
g(x+y)$. The convolution transforms to dot product:  $\widehat{f\ast
  g} = \widehat{f} \cdot \widehat{g}$. The convolution operator is commutative
and associative.

Finally, we need to know Fourier transforms of some simple
functions. The following facts are easily verifiable. Let $f = 1_C$ be
the characteristic function of a linear code $C$. Then $\widehat{f} =
\frac{|C|}{2^n} \cdot 1_{C^{\perp}}$. Let $L(x) = \left\{ \begin{array}{ll}
		    2^n & |x| = 1 \\
		      
 0 &
\mbox{otherwise}
  	              \end{array}
                     \right.$.
Then $\widehat{L}(S) = n - 2|S|$.

Note, for future use, that a code $C \subseteq \H$ has minimal
distance $d$ if and only if $\(1_C \ast
1_C\)(x) = 0$ for $0 < |x| < d$. For a linear code $C$ this is
equivalent to $1_C(x) = 0$  for $0 < |x| < d$.
Observe also that for a function $f$ on the cube holds $Af
= f \ast L$. 

\subsection{The proof of Proposition~\ref{pro:dual}} 
We start with a simple observation that a function supported on a
small set has a large ratio between its second moment and the square
of its first moment. Indeed, let 
a function $F$ be supported on a set $U$. Then, by the Cauchy-Schwarz
inequality, 
\beqn
\label{ineq:c-s}
\E^2 F  = \<F,1_U\>^2 \le \E F^2 \cdot \E \(1_U\)^2 = \frac{|U|}{2^n}
\cdot \E F^2
\eeqn

Hence, to prove (\ref{cover}) it suffices to define a function $F$
supported on $\bigcup_{z \in C'} \(z + B\)$ with $\frac{\E
F^2}{\E^2 F} \le n$. Consider the
adjacency matrix of the subgraph of $\H$ induced by the vertices of $B$. 
Let $f_B$ be an eigenfunction of this matrix corresponding to its maximal
eigenvalue $\l_B$. That is, $f_B$ is supported on $B$ and $\l_B = \frac{\<A
f_B, f_B\>}{\<f_B,f_B\>}$. Since the matrix $A$ is nonnegative, so is
the function $f_B$, and we have $Af_B \ge \l_B f_B$. To see this, note
that $Af_B = \l_B f_B$ on $B$ and, since $Af_B$ is nonnegative, the
inequality holds outside $B$.

For typographic convenience we will write $\l = \l_B$ and $f = f_B$
from now on. 

For a point $z$ in the Hamming cube, let $f_z$ be a shifted version of
$f$, taking $f_z(x) = f(x + z)$. Define 
$$
F = \frac{1}{2^n} \sum_{z \in C'} f_z = 1_{C'} \ast f
$$
This is a nonnegative function supported on $\bigcup_{z \in C'}
\(z + B\)$. We estimate the inner product $\<AF, F\>$ in two ways. 

One one hand, 
$$
AF = F \ast L =  \(1_{C'} \ast f\)
\ast L = 1_{C'} \ast \(f \ast L\) = 
1_{C'} \ast Af \ge \l \(1_{C'} \ast f\) = \l F
$$
Therefore $\<AF, F\> \ge \l_B \<F,F\> = \l_B ~\E F^2$.

On the other hand, by Parseval's identity,
$$
\<AF, F\> = \<\widehat{AF}, \widehat{F}\> = \<\widehat{L} \cdot
\widehat{F},\widehat{F}\> = \<(n-2|S|)\widehat{F},\widehat{F}\> =
\sum_S (n-2|S|)\widehat{F}^2(S)   
$$
Now,  $\widehat{F}(S)
= \widehat{1_{C'}}(S) \cdot \widehat{f}(S) = 0$, for $0 <
|S| < d$. Hence,
$$
\sum_S (n-2|S|)\widehat{F}^2(S) = n \widehat{F}^2(0) + \sum_{|S| \ge
  d} (n-2|S|)\widehat{F}^2(S) \le n \widehat{F}^2(0) + (n-2d)
\sum_{|S| \ge d} \widehat{F}^2(S) \le 
$$
$$
n \widehat{F}^2(0) + (n-2d) 
\sum_S \widehat{F}^2(S) = n \E^2 F + (n-2d) \E F^2
$$
Combining the two estimates on $\<AF,F\>$ and recalling $\l \ge
n-2d+1$, we get
$$
n \E^2 F \ge \(\l - (n-2d)\) \E F^2 \ge \E F^2,
$$
completing the proof. \eprf

\subsection{The proof of Proposition~\ref{pro:ineq}}
\label{subsec:gencod}
The outline of the following proof is
very similar to that of Proposition~\ref{pro:cover}. We suggest that
it is worthwhile to view this proof as a functional version of the
preceding proof. In particular, in light of (\ref{ineq:c-s}) and
(\ref{ineq:es-sup}) below, 
it is useful to define the ``essential support size'' of a function $g$ 
by $2^n \cdot \frac{\E^2 g}{\E g^2}$.

Let $\phi$ be a function on the Hamming cube such that
$\widehat{\phi}^2 = 1_C \ast 1_C$. In other words, $\widehat{\phi \ast
  \phi} = 1_C \ast 1_C$. Since the Fourier
transform on the cube is an involution, up to normalization, we have 
$\phi \ast \phi = 2^n ~\widehat{1_C \ast 1_C} = 2^n
~\widehat{1_C}^2$. What is important is that 
$$
\phi \ast \phi \ge 0~~~ \mbox{and}~~~\frac{\E \phi^2}{\E^2 \phi} =
\frac{\(\phi \ast \phi\)(0)}{\widehat{\phi \ast \phi}(0)} = |C|
$$  
That is, the essential support size of $\phi$ is $\frac{2^n}{|C|}$. 
Note that, for a linear code $C$, we can choose $\phi$ to be (a
multiple of) $1_{C^{\perp}}$.

Take $F = \phi \ast f$. We will show that $\E F^2 \le ~n \E^2
F$. It will take an easy additional step to deduce the desired
inequality $|C| \le n |B|$.

As before, we estimate the inner product $\<AF,F\>$ in two ways.
On one hand,
$$
\<AF,F\> = \<(\phi \ast f) \ast L,\phi \ast f\> =
\<\phi\ast \phi \ast f, f \ast  L\> = 
\<\phi\ast \phi \ast f, A f\> \ge 
$$
$$
\l \<\phi\ast \phi \ast f, f\> = 
\l  \<\phi \ast f, \phi \ast f\> = \l\<F,F\> = \l ~\E F^2
$$
On the other hand, $\<AF,F\> \le n \E^2 F + (n-2d) \E F^2$. The proof
of this fact is exactly the same as in the proof of
Proposition~\ref{pro:cover}, and we omit it.

Combining the two estimates and the assumption $\l \ge n-2d+1$
implies $\E F^2 \le n \E^2 F$. 

Now, $\E^2 F = \E^2 \(\phi \ast f\) = \E^2 \phi ~\E^2 f$. On the other hand, 
$$
\E F^2 = \<F,F\> = \<\phi \ast f, \phi \ast f\> = \<\phi \ast
\phi, f \ast f\> \ge \frac{1}{2^n} \(\phi \ast \phi\)(0) ~\(f
\ast f\)(0) = \frac{1}{2^n} \E \phi^2 ~\E f^2
$$
The inequality follows from nonnegativity of $\phi \ast \phi$.
Since $f$ is supported on $B$, the calculation in (\ref{ineq:c-s}) implies 
\beqn
\label{ineq:es-sup}
|B| \ge 2^n \frac{\E^2 f}{\E f^2} \ge \frac{1}{n} \frac{\E
  \phi^2}{\E^2 \phi} = \frac{1}{n} |C|,
\eeqn
completing the proof. \eprf

\subsection{The proof of Lemma~\ref{lem:eigen-ball}}
We prove the lemma by constructing an explicit function $f$ supported
on $B = B(r)$ with $\frac{\<Af,f\>}{\<f,f\>} \ge \l = 2\sqrt{r(n-r)} -
o(n)$. In fact, we will guarantee more, namely $f \ge 0$ and $Af \ge
\l f$. 

The Hamming ball $B$ contains all the points $x$ of the
Hamming cube with Hamming weight $|x| \le r$. The function $f$ will be
{\it symmetric}, namely its value at a point will depend only on the
Hamming weight of the point. Such a function, of course, is fully
defined by its values $f(0),...,f(n)$ at Hamming
weights $0...n$.  

For a symmetric function $g$ on $\H$ holds
$
Ag(i) = i g(i-1) + (n-i) g(i+1).
$
We start with a preliminary construction of a symmetric function $g$,
setting $g(0) = 1$ and defining $g(i)$ for $1 \le i \le n$ so that the
relation 
\beqn
\label{def:f}
\l g(i) = i g(i-1) + (n-i) g(i+1)
\eeqn
is satisfied
for $i = 1,...,n-1$. We will show below that there exists an integer
$p \le r$ and a real number $\l$ such that 
the function $g$ is nonnegative on the integers $i =
0,...,p$ and nonpositive on $p+1$, and that  $\l \ge 2\sqrt{r(n-r)}
- o(n)$.

Given this, we define $f = g$ for $i = 0,...,p$ and $f = 0$
otherwise. Clearly, $f$ is nonnegative and supported on $B$. 
We claim $Af \ge \l f$. Indeed, by 
definition, $Af(i)  = \l f(i)$, for   
$i \le p - 1$ and for $i > p + 1$. It remains to check the two
boundary values. For  $i = p + 1$,
$Af(i) \ge 0 = \l f(i)$. For $i = p$,
$$
Af(i)  = p f(p-1) + (n-p) f(p+1) = p f(p-1) = p g(p-1) \ge
$$
$$
p g(p-1) + (n-p) g(p+1) =  \l g (p) = \l f(p).
$$
The inequality holds since $g(p+1) \le 0$.

It remains to show $\l \ge 2\sqrt{r(n-r)} - o(n)$. 
We will show that there is a function $r(\l) = (1+ o(n)) \cdot\frac{n -
  \sqrt{n^2-\l^2}}{2} $ such that $g = g_{\l}$ is negative at an
integer point $p \le r(\l)$. Writing $\l$ as a function of
$r$ gives the relation we need, that is $\l \ge 2\sqrt{r(n-r)} - o(n)$.  

Fix $\e > 0$. 
Let $t = \frac{n - \sqrt{n^2-\l^2}}{2}$. We will assume that $g$ is
positive on the interval $\left[0,(1+ \e) t\right]$ and obtain a
contradiction, for a sufficiently large $n$.

By the definition of $g$, 
$$
g(i+1) = \frac{\l g(i) - ig(i-1)}{n-i}.
$$
Set $\theta(i) = \frac{f(i)}{f(i-1)}$. Since $f$ is positive on
$\left[0,(1+ \e) t\right]$ ,
for any $i$ in $\left[t,(1+ \e) t\right]$ holds $\theta(i+1) \ge
\frac{i}{\l} \ge  \frac{t}{\l}$. 
On the other hand, we claim that for any $i > (1+\epsilon/2)t$
holds $\theta(i+1) < \theta(i) \cdot (1-\delta)$, for a positive 
constant $\delta$ depending on $\epsilon$. These two facts evidently
cannot coexist, giving the desired contradiction. 

Indeed, for any $i\ge 0$ holds 
$\theta(i+1) = \frac{\l}{n-i} - \frac{i}{(n-i)\theta(i)}$. We
claim that for 
$i > (1+\epsilon/2)t$, and for any $x > 0$ holds 
$$
\frac{\l}{n-i} - \frac{i}{(n-i)x} < (1 - \delta) x,
$$ 
for some $\delta = \delta(\e) > 0$. In fact, the discriminant of this
quadratic inequality in $x$ is easily seen to be negative, for a
sufficiently small $\delta$. \eprf

\section{Acknowledgements}
We are grateful to Nati Linial for many useful discussions. We also thank
Simon Litsyn and Madhu Sudan for valuable remarks.

\end{document}